\documentclass[12pt]{amsart}

\usepackage[latin1]{inputenc}
\usepackage{natbib}
\usepackage{amssymb}
\usepackage{latexsym}
\usepackage{amsfonts}
\usepackage{amsmath}
\usepackage{bm}
\usepackage{bbm}
\usepackage{mathrsfs}
\usepackage[T1]{fontenc}
\usepackage{ae}
\usepackage[english]{varioref}
\usepackage{graphicx}
\usepackage{a4wide}
\usepackage{enumitem}
\usepackage{extsizes}
\usepackage[usenames, dvipsnames]{color}

\newcommand{\N}{\mathbbm{N}}

\newcommand{\Poi}{\operatorname{Poi}}
\usepackage{amsthm}

\renewcommand{\phi}{\varphi}
\renewcommand{\epsilon}{\varepsilon}

\newcommand\convv{\stackrel{v}{\longrightarrow}}
\newcommand\convd{\stackrel{d}{\longrightarrow}}

\newtheorem{theorem}{Theorem}[section]

\newtheorem{cor}[theorem]{Corollary}

\newtheorem{rmk}[theorem]{Remark}

\begin{document}
%

\numberwithin{equation}{section}
\allowdisplaybreaks

   \author{Anish Ghosh}

   \address{Anish Ghosh, School of Mathematics, Tata Institute of Fundamental Research, Mumbai 400005, India}


   \email{ghosh.anish@gmail.com}


   \author{Maxim S{\o}lund Kirsebom}
   \address{Maxim S{\o}lund Kirsebom, Department of Mathematics, University of Hamburg, 20146 Hamburg, Germany}
   \email{maximkirsebom@gmail.com}


    \author{Parthanil Roy}

   \address{Parthanil Roy, Theoretical Statistics and Mathematics Unit, Indian Statistical Institute, Bangalore 560059, India}


   \email{parthanil.roy@gmail.com}


\title[Continued fractions and extreme value theory]{Continued fractions, the Chen-Stein method and extreme value theory}

\thanks{A.\ G.\ gratefully acknowledges support from a grant from the Indo-French Centre for the Promotion of Advanced Research; a Swarnajayanti Fellowship from Department of Science and Technology, Government of India and a MATRICS grant from the Science and Engineering Research Board. P.\ R.\ acknowledges the support from a MATRICS grant from the Science and Engineering Research Board and a Swarnajayanti Fellowship from Department of Science and Technology, Government of India.}

\begin{abstract}
In this work, we deal with extreme value theory in the context of continued fractions using techniques from probability theory, ergodic theory and real analysis. We give an upper bound for the rate of convergence in the Doeblin-Iosifescu asymptotics for the exceedances of digits obtained from the regular continued fraction expansion of a number chosen randomly from $(0,1)$ according to the Gauss measure. As a consequence, we significantly improve the best known upper bound on the rate of convergence of the maxima in this case. We observe that the asymptotics of order statistics and the extremal point process can also be investigated using our methods. 
\end{abstract}



\maketitle

\section{introduction}

This short paper establishes an upper bound for the Doeblin-Iosifescu asymptotics for exceedances (defined below) arising from the Gauss dynamical system. We briefly recall the basic facts about continued fraction expansions and the Gauss map. The reader is referred to the classic text \cite{khintchine:1964} for more details. Let $X = (0,1)$ and for all $x \in X,$ let $[A_1(x), A_2(x), \ldots]$ denote the regular continued fraction expansion. Define a transformation $T: X \to X$ by
\begin{equation} \label{def:T}
T(x)=\{1/x\},
\end{equation}
where  $\{\cdot\}$ denotes the fractional part. With the notations above, for all $x \in X\backslash \mathbbm{Q}$,
\begin{align*}
A_1(x) &= 1/x - T(x), \; \\
A_{j+1}(x) &= A_j(T(x)) = A_1(T^j(x)) \mbox{ for all }j \in \mathbb{N}.
\end{align*}
It is easy to check that $T$ defines a nonsingular transformation on $(X, \lambda)$, where $\lambda$ denotes the Lebesgue measure. This means that  for all measurable $B \subseteq X$, we have $\lambda \circ T^{-1}(B) =0$ if and only if $\lambda(B) =0$.

Let $\widehat{T}: L^1(X, \lambda) \to L^1(X, \lambda)$ denote the dual operator (see, for example, Page~33 of \cite{aaronson:1997}) corresponding to $T$ that satisfies
\[
\int_X \widehat{T}(f) g d\lambda = \int_X f (g \circ T) d\lambda
\]
for all $f \in L^1(X, \lambda)$ and for all $g \in L^\infty(X, \lambda)$. It is easy to extend the domain of definition of $\widehat{T}$ to all nonnegative measurable functions. Solving the functional equation $\widehat{T}(h) = h$, we get $h(x) = (1+x)^{-1} \in L^\infty(X, \lambda)$. Hence by Proposition~1.4.1 of \cite{aaronson:1997}, the probability measure $$P(dx) = ((1+x)\log{2})^{-1}dx$$ on $X$ is $T$-invariant making $T$ a positive transformation (see, for example, \cite{aaronson:1997}). The measure $P$ is known as the Gauss measure.

From now on, we shall think of $\{A_n\}_{n \geq 1}$ as a sequence of random variables $A_n:X\to\N$ defined on the probability space $(X, P)$. The $T$-invariance of $P$ makes this a stationary sequence, i.e., for all $k, l \in \mathbb{N}$, for all $m_1, m_2, \ldots, m_k \in \mathbb{N}$ and for all Borel subset $B \subseteq X^k$,
\[
P\big( (A_{m_1}, A_{m_2}, \ldots, A_{m_k}) \in B\big) = P\big( (A_{m_1+l}, A_{m_2+l}, \ldots, A_{m_k+l}) \in B\big).
\]
We are interested in the extreme value theory for this stationary stochastic process.
To the best of our knowledge, the first work in this direction was carried out by \cite{doeblin:1940}, who, among many other results, rightly observed that exceedances have Poissonian asymptotics: for all $u>0$,
\begin{equation}
\mathcal{E}^u_n := \# \{1 \leq i \leq n:\, A_i \log{2} > nu\}  \convd \mathcal{E}^u_\ast \sim \Poi(u^{-1})\label{conv:doeblin_iosifescu}
\end{equation}
under $P$. Here $\convd$ denotes convergence in distribution and the notation $\mathcal{E}^u_\ast \sim \Poi(u^{-1})$ means that
\begin{equation*}
	P\big( \mathcal{E}^u_\ast = k\big)=\frac{u^{-k}e^{-u^{-1}}}{k!},\,\;k =0,1, 2, \ldots\,.
\end{equation*}
 However, Doeblin's proof of \eqref{conv:doeblin_iosifescu} had a subtle error, which was corrected much later in Theorem~2 of \cite{iosifescu:1977}. Therefore, we shall refer to \eqref{conv:doeblin_iosifescu} as the Doeblin-Iosifescu asymptotics; they form the background of this paper.

Seemingly unaware of the work of \cite{doeblin:1940}, three decades later \cite{galambos:1972} showed that for all $u > 0$,
\begin{equation} \label{conv:maxima}
P\left( \frac{\log{2} \, \max_{i=1}^n A_i}{n} \leq u\right) \to e^{-u^{-1}},
\end{equation}
which is a restatement of $P(\mathcal{E}^u_n = 0) \to P(\mathcal{E}^u_\ast = 0)$ and hence an easy consequence of \eqref{conv:doeblin_iosifescu}. However, because of the subtle mistake of \cite{doeblin:1940}, the above result of \cite{galambos:1972} stands as the first correctly proven result on extreme value theory of continued fractions. This has remained a topic of current interest; see, for example, the generalizations of \eqref{conv:maxima} to fibred systems by \cite{nakada:natsui:2003} and to Oppenheim continued fractions by \cite{chang:ma:2017}.

In view of the above, the following question arises naturally:\\

 \noindent \emph{What is the rate of convergence in the of the asymptotics in  \eqref{conv:doeblin_iosifescu}?}\\

 \noindent In this paper, we give an upper bound on the rate of convergence using the Chen-Stein method of \cite{arratia:goldstein:gordon:1989} (more specifically, Theorem~\ref{thm_Arratia_et_al} below). As far as we are aware, our work is the first to specifically employ the Chen-Stein method in the context of Gauss map and continued fractions.

 The Chen-Stein method is a very useful technique which yields an upper bound that is uniform in $u$ bounded away from zero; see, Theorem~\ref{thm:rate_DI} below. As a consequence, we also get a locally uniform (in $(0, \infty]$) upper bound for the convergence of distribution functions in \eqref{conv:maxima} and this bound is much better than the best known bound given in \cite{philipp:1976} (we improve a slowly varying rate of convergence to a polynomial one; see Remark~\ref{rmk:comparison} below). In fact, we give a bound on the rate of convergence of the $k^{th}$ maxima, not just the maxima, and the Chen-Stein method is  powerful enough to ensure that this locally uniform upper bound turns out to be uniform over $k \in \mathbb{N}$ as well (see Corollary~\ref{cor:rate_extremes}).

Note that \eqref{conv:maxima} implies that the $A_i$'s are in the Fr\'{e}chet(1) maximal domain of attraction.  It is not difficult to observe that \eqref{conv:maxima} holds because the $A_i$'s enjoy a very strong exponential mixing property (see \eqref{mixing_exp} below), and each $A_i$ (which are anyway identically distributed because of stationarity) is regularly varying with index $-1$, i.e.,
\begin{equation}
nP\left(\frac{A_1\log{2}}{n} \in \cdot \right) \convv \nu  \label{conv:vague:A_i}
\end{equation}
as measures on $(0, \infty]$. Here ``$\convv$'' denotes vague convergence and $\nu$ is the unique measure on $(0, \infty]$ satisfying $\nu\big((u, \infty]\big) = u^{-1}$ for all $u \in (0, \infty)$. This was essentially the proof given in \cite{galambos:1972} except that he did not use the language of vague convergence, and presented a direct proof instead.

The above vague convergence will play a very important role in this paper. Since $A_1$ is an integer-valued random variable, it follows that for each $u >0$,
\begin{align}
P\left(\frac{A_1\log{2}}{n}>u\right) &= P\left(A_1 \geq \left\lceil \frac{nu}{\log{2}}\right\rceil\right), \nonumber
\intertext{from which, using Lemma~1 of \cite{galambos:1972}, we get}
nP\left(\frac{A_1\log{2}}{n}>u\right) &= \frac{n\log{\left(1+\frac{1}{\lceil n u/\log{2} \rceil}\right)}}{\log{2}} \to u^{-1} \label{conv:upper_tail}
\end{align}
as $n \to \infty$. From the above convergence, \eqref{conv:vague:A_i} follows by invoking Theorem~3.6 of \cite{resnick:2007}. Further, using the inequality $\log{(1+x)} \leq x$ whenever $x>0$, we get the following upper bound, which will also be very useful in this paper: for all $u>0$,
\begin{equation}
P\left(\frac{A_1\log{2}}{n}>u\right) = \frac{\log{\left(1+\frac{1}{\lceil n u/\log{2} \rceil}\right)}}{\log{2}}\leq  \frac{\log{\left(1+\frac{\log{2}}{nu}\right)}}{\log{2}} \leq \frac{1}{nu}. \label{imp_upper_bound}
\end{equation}

In some sense, the $A_i$'s behave very much like an i.i.d.~sequence because of the following exponential mixing property. For all $m, n \in \mathbb{N}$, for all $F \in \sigma(A_1, A_2, \ldots, A_m)$, and for all $H \in \sigma(A_{m+n}, A_{m+n+1}, \ldots)$,
\begin{equation}
|P(F \cap H)-P(F)P(H)| \leq \psi(n) P(F)P(H) = C \theta^{-n}P(F)P(H), \label{mixing_exp}
\end{equation}
where $\psi(n) = C \theta^{-n}$ for some $C>0$ and $\theta >1$; see, for example, Lemma~2 of \cite{galambos:1972}.

In order to state our main result and its corollary, we need to introduce some notation as described below. For each $n \in \mathbb{N}$ and for each $k \in \{1, 2, \ldots, n\}$, denote by $M_n^{(k)}$, the $k^{th}$ largest in the set $\{A_i\log{2}: 1 \leq i \leq n\}$. Then it follows from \eqref{conv:doeblin_iosifescu} that for all $u>0$,
\[
P\left(\frac{M_n^{(k)}}{n} \leq u\right) = P(\mathcal{E}^u_n \leq k-1) \rightarrow P(\mathcal{E}^u_\ast \leq k-1) = e^{-u^{-1}}\sum_{i=0}^{k-1} \frac{u^{-i}}{i!}.
\]
Obviously, the $k=1$ case has already been taken care of in \eqref{conv:maxima} above. Also, let $\{l_n\}$ be a sequence of positive real numbers such that
\begin{equation}
l_n \theta^{l_n} = n \label{defn_of_l_n}
\end{equation}
for all $n \in \mathbb{N}$ (here $\theta$ is as in \eqref{mixing_exp} above). Clearly, such a sequence exists by the intermediate value theorem and it increases to infinity at a rate strictly slower than $\log{n}$.

We are now ready to state our main result.

\begin{theorem}\label{thm:rate_DI}
With the notation as above, we have the following upper bound on the rate of convergence in \eqref{conv:doeblin_iosifescu}: there exists $\kappa >0$ such that for all $\delta >0$ and for all $n \in \mathbb{N}$,
\begin{equation}
\sup_{u \in [\delta, \infty)} d_{TV}(\mathcal{E}^u_n, \mathcal{E}^u_\ast) \leq  \,\frac{\kappa}{\min{\{\delta, \delta^2\}}}\,\frac{l_n}{n}, \nonumber
\end{equation}
where $d_{TV}$ denotes the total variation distance.
\end{theorem}

\noindent We would like to mention that we blend probability theory (namely, the Chen-Stein method; see Theorem~\ref{thm_Arratia_et_al}), ergodic theory (specifically, the exponential mixing property \eqref{mixing_exp}) and real analysis (more precisely, a second order regular variation estimate; the second inequality in \eqref{est:RV2}) to prove the result above. 

Theorem~\ref{thm:rate_DI} has the following very strong consequence on the rate of weak convergence of scaled $k^{th}$ maxima. The upper bound here is uniform over $u$ bounded away from zero and uniform over $k \in \mathbb{N}$ at the same time.

\begin{cor} \label{cor:rate_extremes} With $\kappa$ as in Theorem~\ref{thm:rate_DI}, we get that for all $\delta>0$ and for all $n \in \mathbb{N}$,
\begin{equation}
\sup_{k \in \mathbb{N}}\sup_{u \in [\delta, \infty)} \left|P\left(\frac{M_n^{(k)}}{n} \leq u\right) - e^{-u^{-1}}\sum_{i=0}^{k-1} \frac{u^{-i}}{i!}\right| \leq  \,\frac{\kappa}{\min{\{\delta, \delta^2\}}}\,\frac{l_n}{n}. \nonumber
\end{equation}
\end{cor}
\noindent The above corollary follows from Theorem~\ref{thm:rate_DI} by restricting the supremum in the definition of total variation distance to sets of the form $\{0, 1, \ldots, k-1\}$ with $k$ running over the set of all positive integers.

\begin{rmk}\label{rmk:iid} \textnormal{Note that if $A_i$'s were i.i.d.\ with same marginal distribution, then by \cite{resnick:dehaan:1989}, we would have obtained an upper bound of $O\left(\frac{1}{n}\right)$ on the rate of convergence of the maxima sequence. The Chen-Stein method gives the same rate in the i.i.d.\ case. In the Gauss dynamical system, we get an extra factor of $l_n$ because of the dependence of the $A_i$'s. However, since $l_n = o(\log{n})$, it follows that our bound on the rate of convergence is $o\left(\frac{\log{n}}{n}\right)$. Therefore, we almost attain the rate obtained in the i.i.d.\ case.}
\end{rmk}

\begin{rmk}\label{rmk:comparison}
\noindent \textnormal{The best known rate of convergence for the maxima in our setup was obtained by \cite{philipp:1976}, who gave an upper bound of $O\left(e^{-(\log{n})^\delta}\right)$ with $\delta \in (0,1)$ (the constant in $O$ depends on $\delta$). Note that $e^{-(\log{n})^\delta}$ is a slowly varying function of $n$. Therefore, by the Potter bound (see, for example, Page~32 of \cite{resnick:2007}), it follows that
$
n^{-\eta} =o\left(e^{-(\log{n})^\delta}\right)
$
for all $\eta >0$ and for all $\delta \in (0,1)$. Hence, by Remark~\ref{rmk:iid}, it follows that
\[
\frac{l_n}{n} = o\left(e^{-(\log{n})^\delta}\right)
\]
for all $\delta \in (0,1)$. Therefore, our bound on the rate of convergence is significantly better than the one obtained by \cite{philipp:1976}. More precisely, we improve a slowly varying rate of convergence to a polynomial one, bettering an error term that was used by Philipp in his proof of a conjecture of Paul Erd\"{o}s.}
\end{rmk}

Note that the $D$ and $D^\prime$ conditions of \cite{davis:1983} follow from \eqref{mixing_exp}. Therefore, by Example~5.1 in \cite{davis:hsing:1995}, the following extremal point process weak convergence holds in the space $\mathcal{M}_p((0, \infty])$ of all Radon point measures (on $(0, \infty]$) equipped with the vague metric:
\begin{equation}
  Q_n := \sum_{i=1}^n \delta_{\frac{A_i \log{2}}{n}} \convd Q_\ast \sim PRM((0, \infty], \nu).  \label{conv:weak:Q_n:PRM}
\end{equation}
Here the limit $Q_\ast$ is a Poisson random measure on $(0, \infty]$ with mean measure $\nu$; see Section~4.1 of \cite{tyran-kaminska:2010} for a direct proof of \eqref{conv:weak:Q_n:PRM}. In this paper, we observe that a tiny detour of our proof of Theorem~\ref{thm:rate_DI} yields \eqref{conv:weak:Q_n:PRM}; see Section~\ref{sec:ppconv} below.

\subsection*{Acknowledgements} This work was initiated during a visit by M.K. and P.R. at the \emph{Tata Institute of Fundamental Research, Mumbai} and a significant portion of the work was carried out when the authors were at the \emph{International Centre for Theoretical Sciences, Bangalore} for the program \emph{Probabilistic Methods in Negative Curvature} (ICTS/pmnc2019/03). We thank both institutes for their hospitality and the lovely working conditions. We would also like to acknowledge an anonymous reviewer and an executive editor for their careful reading of the paper. Their detailed comments have significantly improved our work (especially, Remarks~\ref{rmk:iid} and \ref{rmk:rs}). 

\section{Proofs} \label{sec:proofs}

As mentioned earlier, the proof of Theorem~\ref{thm:rate_DI} relies on the Chen-Stein method of \cite{arratia:goldstein:gordon:1989}. We first state their result and then present our proof. Finally, we observe how a tiny detour of the proof also establishes the weak convergence of the extremal point process of the digits arising in the continued fraction expansion.

\subsection{The Chen-Stein Method of \cite{arratia:goldstein:gordon:1989}}

Let $\mathcal{I}$ be an index set and $\{X_\alpha \sim Ber(p_\alpha)\}_{\alpha \in \mathcal{I}}$ be a collection of possibly dependent Bernoulli random variables. Suppose, for each $\alpha \in \mathcal{I}$, there exists a subset $B_\alpha \subseteq \mathcal{I}$ such that roughly speaking, $X_\alpha$ is nearly independent of $\{X_\beta:\, \beta \in \mathcal{I} \setminus B_\alpha\}$. \cite{arratia:goldstein:gordon:1989} called $B_\alpha$ the ``neighborhood of dependence'' of $X_\alpha$. Following their notation, we define
\begin{align}
b_1  &=  \sum_{\alpha \in \mathcal{I}} \sum_{\beta \in B_\alpha} p_\alpha p_\beta,  \label{defn:b_1}\\
b_2 &=  \sum_{\alpha \in \mathcal{I}} \sum_{\beta \in B_\alpha \setminus \{\alpha\}} E(X_\alpha X_\beta), \label{defn:b_2}\\
b_3 &=  \sum_{\alpha \in \mathcal{I}} E\big[| E(X_\alpha - p_\alpha | \mathcal{H}_\alpha)|\big], \label{defn:b_3}
\end{align}
where $\mathcal{H}_\alpha$ is the $\sigma$-field generated by $\{X_\beta:\, \beta \in \mathcal{I} \setminus B_\alpha\}$.

\begin{theorem}[Theorem~2 of \cite{arratia:goldstein:gordon:1989}] \label{thm_Arratia_et_al} Partition $\mathcal{I}$ into disjoint nonempty subsets $\mathcal{I}_1, \mathcal{I}_2, \ldots , \mathcal{I}_k$. Let $\{Y_\alpha \sim Poi(p_\alpha)\}_{\alpha \in \mathcal{I}}$ be a collection of independent Poisson random variables. Set
\[
W_j = \sum_{\alpha \in \mathcal{I}_j} X_\alpha, \;\;\; Z_j = \sum_{\alpha \in \mathcal{I}_j} Y_\alpha \;\;\; \mbox{and}  \;\;\;\lambda_j = \sum_{\alpha \in \mathcal{I}_j} p_\alpha.
\]
 Then
 \begin{align}
 & d_{TV}\big(\mathcal{L}(W_1, W_2, \ldots, W_k), \mathcal{L}(Z_1, Z_2, \ldots, Z_k)\big) \nonumber\\
 &\;\;\;\;\; \;\;\;\;\; \;\;\;\;\; \leq \min{\left\{2, 2.8\max_{1 \leq j \leq k} \lambda_j^{-1/2}\right\}}(2b_1 + 2b_2 + b_3),  \label{bound_Arratia_et_al}
 \end{align}
where $\mathcal{L}$ denotes the joint law.
\end{theorem}

We would like to elaborate a bit on the phrase ``nearly independent'' used above in the context of neighborhood of dependence $B_\alpha$. In many examples (e.g., $m$-dependent time-series models, certain random graph asymptotics, etc.) where Theorem~\ref{thm_Arratia_et_al} is
used, $X_\alpha$ is totally independent of $\{X_\beta: \beta \not\in B_\alpha\}$ making $b_3 = 0$. In our case, however, we need to bound $b_3$ tightly using the ``near independence'' property \eqref{mixing_exp}.

\subsection{Proof of Theorem~\ref{thm:rate_DI}} \label{sec:mainproof}

Define a new Poisson random variable $\tilde{\mathcal{E}}^u_n$ with mean $n P(n^{-1}A_1\log{2} > u)$. The basic strategy of the proof is to use that
\begin{equation}\label{d_TV_triangle_ineq}
	d_{TV}(\mathcal{E}^u_n, \mathcal{E}^u_\ast)\leq d_{TV}(\mathcal{E}^u_n, \tilde{\mathcal{E}}^u_n)+ d_{TV}(\tilde{\mathcal{E}}^u_n, \mathcal{E}^u_\ast)
\end{equation}
and to estimate each term separately. The bound on $d_{TV}(\mathcal{E}^u_n, \tilde{\mathcal{E}}^u_n)$ will need Chen-Stein method and the exponential mixing property \eqref{mixing_exp} while the second term $d_{TV}(\tilde{\mathcal{E}}^u_n, \mathcal{E}^u_\ast)$ will be estimated using a hard analytic bound on the second order term of the convergence in \eqref{conv:upper_tail}. Thus, our proof combines tools from probability theory, ergodic theory and real analysis in a systematic manner.

We will first show that there exists $\kappa_1>0$ such that for all $u>0$ and for all $n \geq 1$,
\begin{equation}
d_{TV}(\mathcal{E}^u_n, \tilde{\mathcal{E}}^u_n) \leq \kappa_1 \max\left\{\frac{1}{u}, \frac{1}{u^2}\right\} \frac{l_n}{n}, \label{conv:DI_step_1}
\end{equation}
where $l_n$ is as in \eqref{defn_of_l_n}.
To this end, set
\begin{equation}
D= (u, \infty]. \label{defn:D_exceed}
\end{equation}
We shall use Theorem~\ref{thm_Arratia_et_al} with $\mathcal{I}= \{1, 2, \ldots, n\}$, $k=1$, $X_\alpha = I_{(n^{-1}A_\alpha \log{2} \in D)}$ (and hence $p_\alpha = E(X_\alpha) = E(X_1) = P(n^{-1}A_1 \log{2} \in D)$) and $B_\alpha = (\alpha -l_n, \alpha + l_n) \cap \mathcal{I}$
for each $\alpha \in \mathcal{I}$. Note that with these choices we have $W_1=\mathcal{E}^u_n$ and $Z_1$ may be thought of, intuitively, as ``$W_1$ if the $X_{\alpha}$'s were independent''.

Because of stationarity, we get
$$
\sum_{\alpha \in \mathcal{I}} p_\alpha = n P(n^{-1}A_1 \log{2} \in D) = n P(n^{-1}A_1\log{2} > u).
$$
In order to establish \eqref{conv:DI_step_1}, we have to estimate the quantities defined by \eqref{defn:b_1}, \eqref{defn:b_2} and \eqref{defn:b_3}. For the first one, observe that
\begin{align}
b_1
&=  \sum_{\alpha \in \mathcal{I}} \sum_{\beta \in B_\alpha} p_\alpha p_\beta \nonumber\\
&=  \sum_{\alpha=1}^n \sum_{\beta \in B_\alpha} \big(P(n^{-1}A_1 \log{2} \in D)\big)^2, \nonumber\\
\intertext{from which, using \eqref{imp_upper_bound}, we get}
b_1& \leq 2n l_n \left(\frac{1}{nu}\right)^2 = \frac{2}{u^2} \frac{l_n}{n}. \label{bound_on_b1}
\end{align}

In order to bound the second term in \eqref{bound_Arratia_et_al}, note that for any $\alpha, \beta \in \mathbb{N}$ such that $\alpha \neq \beta$,
\begin{align*}
E(X_\alpha X_\beta)
&= P\big(n^{-1}A_\alpha \log{2} \in D, n^{-1}A_\beta \log{2} \in D\big) \\
& \leq (1+\psi(|\alpha - \beta|)) P\big(n^{-1}A_\alpha \log{2} \in D)P(n^{-1}A_\beta \log{2} \in D\big),
\end{align*}
where the last step follows from \eqref{mixing_exp}. Applying stationarity, \eqref{imp_upper_bound} and the inequality $\psi(n) \leq C$, we get from the above bound that
\[
E(X_\alpha X_\beta)
\leq (1+C) \big(P\big(n^{-1}A_1 \log{2} \in D)\big)^2 \leq (1+C) \left(\frac{1}{nu}\right)^2
\]
for all $\alpha \neq \beta$. Hence
\begin{equation}
b_2 = \sum_{\alpha \in \mathcal{I}} \sum_{\beta \in B_\alpha \setminus \{\alpha\}} E(X_\alpha X_\beta) \leq 2n l_n (1+C)  \frac{1}{u^2n^2} = \frac{2(1+C)}{u^2} \frac{l_n}{n}. \label{bound_on_b2}
\end{equation}

Finally, we need to estimate \eqref{defn:b_3}. Fixing $\alpha \in \mathcal{I}$ and taking $F=(n^{-1}A_\alpha \log{2} \in D)$ with $D$ as in \eqref{defn:D_exceed}, we see that \eqref{mixing_exp} yields
\[
p_\alpha P(H)(1-\psi(l_n)) \leq
P\big((n^{-1}A_\alpha \log{2} \in D) \cap H\big) \leq p_\alpha P(H)(1+\psi(l_n))
\]
for all $H \in \mathcal{H}_\alpha=\sigma\{X_\beta:\, \beta \in \mathcal{I} \setminus B_\alpha\}$. The above pair of inequalities can be rewritten as
\[
p_\alpha P(H)(1-\psi(l_n)) \leq
\int_H X_\alpha dP = \int_H E(X_\alpha|\mathcal{H}_\alpha) dP \leq p_\alpha P(H)(1+\psi(l_n))
\]
yielding
\[
-\int_H p_\alpha \, \psi(l_n)\,dP \leq
\int_H E(X_\alpha- p_\alpha |\mathcal{H}_\alpha)\,dP \leq \int_H p_\alpha \, \psi(l_n)\,dP,
\]
which holds for all $H \in \mathcal{H}_\alpha$ and hence
\[
\big|E(X_\alpha- p_\alpha |\mathcal{H}_\alpha)\big| \leq p_\alpha \, \psi(l_n)= P(n^{-1}A_1 \log{2} \in D) \,\psi(l_n)
\]
almost surely. Therefore, we get
\begin{align*}
b_3 &=  \sum_{\alpha \in \mathcal{I}} E\big[| E(X_\alpha - p_\alpha | \mathcal{H}_\alpha)|\big]\\
       & \leq n P(n^{-1}A_1 \log{2} \in D)\psi(l_n) \\
       & \leq n \frac{1}{nu} C \theta^{-l_n}\\
       & = \frac{C }{u}\frac{l_n}{n},
\end{align*}
where we used \eqref{imp_upper_bound} and the last step follows from the choice of $l_n$ as given in \eqref{defn_of_l_n}. The above upper bound, along with \eqref{bound_on_b1} and \eqref{bound_on_b2}, yields \eqref{conv:DI_step_1} thanks to Theorem~\ref{thm_Arratia_et_al}.

We now move on to estimating the second term in \eqref{d_TV_triangle_ineq}. We first use Taylor's theorem to obtain the inequality $|\log(1+x)-x| \leq \frac{x^2}{2}$, which can be rewritten as
\begin{equation}
\left|\frac{\log{(1+x)}}{x}-1\right| \leq \frac{x}{2} \label{imp_upper_bound_2}
\end{equation}
for all $x>0$. Using this inequality, we shall now bound the second order term of the convergence in \eqref{conv:upper_tail}.

To this end, note that
\begin{align}
&|n P(n^{-1} A_1 \log{2} \in D) - u^{-1}| = \left|\frac{n \log{(1+\lceil n u/\log{2}\rceil^{-1})}}{\log{2}} - \frac{1}{u} \right| \nonumber\\
& \leq \frac{n}{\log{2}}\left|\log{\left(1+\frac{1}{\lceil n u/\log{2}\rceil}\right)} - \log{\left(1+\frac{1}{n u/\log{2}}\right)} \right| \nonumber\\
&\hspace{4.5cm} + \; \; \; \frac{1}{u} \left|\frac{\log{(1+n^{-1}u^{-1}\log{2})}}{n^{-1}u^{-1}\log{2}} -1 \right|.\nonumber
\end{align}
By virtue of \eqref{imp_upper_bound_2}, the second term above is bounded by
$
\frac{\log{2}}{2u^2n}.
$
On the other hand, using the mean value theorem, we can estimate the first term as follows:
\begin{eqnarray*}
\frac{n}{\log{2}}\left|\log{\left(1+\frac{1}{\lceil n u/\log{2}\rceil}\right)} - \log{\left(1+\frac{1}{n u/\log{2}}\right)} \right| \leq \frac{n}{\log 2} \frac{{(\log 2)}^2}{n^2 u^2}=\frac{\log 2}{u^2n}
\end{eqnarray*}

Therefore, by Lemma~(8) of \cite{freedman:1974}, it follows that
\begin{equation}
d_{TV}(\tilde{\mathcal{E}}^u_n, \mathcal{E}^u_\ast) \leq |n P(n^{-1} A_1 \log{2} \in D) - u^{-1}| \leq \frac{3\log{2}}{2u^2} \frac{1}{n}. \label{est:RV2}
\end{equation}
The above inequality, \eqref{conv:DI_step_1} and \eqref{d_TV_triangle_ineq} imply that there exists a constant $\kappa \in (0, \infty)$ such that for all $u>0$ and for all $n \geq 1$,
\[
d_{TV}(\mathcal{E}^u_n, \mathcal{E}^u_\ast) \leq \kappa \max\left\{\frac{1}{u}, \frac{1}{u^2}\right\} \frac{l_n}{n},
\]
from which Theorem~\ref{thm:rate_DI} follows. 

\begin{rmk}\label{rmk:rs} \textnormal{We would like to mention here an alternative approach pointed out to us by an anonymous referee. Namely, Theorem~1 of \cite{smith:1988} gives a similar Chen-Stein type upper bound in the more general setup of non-stationary processes. It is possible to use this result to give a bound on $d_{TV}(\mathcal{E}^u_n, \tilde{\mathcal{E}}^u_n)$ in our work leaving the estimation of $d_{TV}(\tilde{\mathcal{E}}^u_n, \mathcal{E}^u_\ast)$ (based on hard analysis) as it is. This will involve (in the notation of \cite{smith:1988}) coming up with the function $g(n, r)$, the subsets $I_{nrk}$ and $I^{\ast}_{nrk}~(\subseteq I_{nrk}$), the latter being very similar to a neighborhood of dependence, and verifying the Condition $D^\prime$ of \cite{smith:1988}. We think that this will be more involved than the estimation of the terms $b_1$ and $b_2$ of our paper. On the other hand, Condition $D$ of \cite{smith:1988} will follow directly from the exponential mixing property~(1.7) of our paper and this verification will be shorter than bounding the term $b_3$ in our work. Overall, we feel that an application of Theorem~1 of \cite{smith:1988}, instead of Theorem~2 of \cite{arratia:goldstein:gordon:1989}, will perhaps result in an argument of similar length. However, we have not compared the rates obtained by these two results in our setup.}
\end{rmk}

\subsection{New Proof of \eqref{conv:weak:Q_n:PRM}} \label{sec:ppconv}

By Theorem~4.7 of \cite{kallenberg:1983}, in order to establish \eqref{conv:weak:Q_n:PRM}, it is enough to show the following:
\begin{enumerate}
  \item[(i)] For all $u, v \in (0, \infty]$ with $u < v$,
  \begin{equation}
  E\big(Q_n((u, v])\big) \longrightarrow E\big(Q_\ast((u, v])\big)= \nu((u, v]) = u^{-1} - v^{-1}\label{conv:weak:Q_n:PRM_step_1}
  \end{equation}
  as $n \to \infty$. Of course, we follow the convention $\infty^{-1}=0$.

  \item[(ii)] Whenever $ 0 < u_1 < v_1 < u_2 < v_2 < \cdots < u_k < v_k \leq \infty$,
    \begin{align}
  &P\big(Q_n((u_1, v_1])=0, Q_n((u_2, v_2])=0, \ldots, Q_n((u_k, v_k])=0\big) \nonumber\\
  &\; \; \rightarrow P\big(Q_\ast((u_1, v_1])=0, Q_\ast((u_2, v_2])=0, \ldots, Q_\ast((u_k, v_k])=0\big) \label{conv:weak:Q_n:PRM_step_2}\\
  &\; \; \;\;\; = \prod_{i=1}^{k} e^{-\nu((u_i, v_i])} = \exp{\left\{-\nu\left(\bigcup_{i=1}^{k} (u_i, v_i]\right)\right\}} \nonumber
   \end{align}
as $n \to \infty$.
\end{enumerate}

By linearity of expectation, in order to establish \eqref{conv:weak:Q_n:PRM_step_1}, it is enough to do so with $u \in (0, \infty)$ and $v =\infty$. This special case follows using stationarity of $A_i$'s and \eqref{conv:upper_tail} as shown below:
\begin{align*}
E\big(Q_n((u, \infty])\big)
&= E\left[\sum_{i=1}^n \delta_{n^{-1}A_i \log{2}}((u, \infty])\right]\\
&=n P\left(n^{-1}A_1 \log{2}>u\right) \\
&\to \nu((u,\infty]),
\end{align*}
as $n \to \infty$). This proves \eqref{conv:weak:Q_n:PRM_step_1}.

On the other hand, verification of \eqref{conv:weak:Q_n:PRM_step_2} will need a tiny detour of the proof of Theorem~\ref{thm:rate_DI} (as carried out in \cite{chiarini:cipriani:hazra:2015} in the context of Gaussian free fields) and Theorem~\ref{thm_Arratia_et_al} will again play a significant role in the proof. To this end, fix $0 < u_1 < v_1 < u_2 < v_2 < \cdots < u_k < v_k \leq \infty$ and set
\begin{equation}
D= \bigcup_{i=1}^k \,(u_i, v_i]. \label{defn:D_2}
\end{equation}
Note that by \eqref{conv:vague:A_i} and Proposition~3.12 of \cite{resnick:1987}, it follows that
\begin{equation}
nP\left(n^{-1}A_1 \log{2} \in D\right) \to \nu(D) \label{imp_impli_vague_conv}
\end{equation}
as $n \to \infty$ for $D$ as in \eqref{defn:D_2}. Therefore by changing the definition of $D$ from \eqref{defn:D_exceed} to \eqref{defn:D_2} in the proof of Theorem~\ref{thm:rate_DI} and using \eqref{imp_impli_vague_conv}, it is easy to show that
\[
d_{TV}(Q_n(D), Q_\ast(D)) \to 0
\]
as $n \to \infty$. In particular, $P(Q_n(D)=0) \to P(Q_\ast(D)=0) = e^{-\nu(D)}$, which is a restatement of \eqref{conv:weak:Q_n:PRM_step_2}. This completes the proof of \eqref{conv:weak:Q_n:PRM} based on the Chen-Stein method of \cite{arratia:goldstein:gordon:1989}.


\begin{thebibliography}{9}
\expandafter\ifx\csname natexlab\endcsname\relax\def\natexlab#1{#1}\fi

\bibitem[Aaronson(1997)]{aaronson:1997}
{\sc J.~Aaronson} (1997): {\em An Introduction to Infinite Ergodic Theory\/}.
\newblock American Mathematical Society, Providence.

\bibitem[Arratia et~al.(1989)]{arratia:goldstein:gordon:1989}
{\sc R.~Arratia, L.~Goldstein {\rm and} L.~Gordon} (1989): Two moments suffice
  for {P}oisson approximations: the {C}hen-{S}tein method.
\newblock {\em Ann. Probab.\/} 17:9--25.

\bibitem[Chang and Ma(2017)]{chang:ma:2017} {\sc Y.~Chang and J.~Ma} (2017): Some distribution results of the Oppenheim continued
fractions.
\newblock{Monatsh. Math.\/} 184(3): 379--399.

\bibitem[Chiarini et~al.(2015)]{chiarini:cipriani:hazra:2015}
{\sc A.~Chiarini, A.~Cipriani {\rm and} R.~S. Hazra} (2015): A note on the extremal process of the supercritical {G}aussian free field.
\newblock {\em Electron. Commun. Probab.\/} 20:paper no. 74, 10 pages.

\bibitem[Davis(1983)]{davis:1983}
{\sc R.~Davis} (1983): Stable limits for partial sums of dependent random variables.
\newblock {\em Ann. Probab.\/} 11(2):262--269.

\bibitem[Davis and Hsing(1995)]{davis:hsing:1995}
{\sc R.~Davis {\rm and} T.~Hsing} (1995): Point processes for partial sum convergence for weakly dependent random variables with infinite variance.
\newblock {\em Ann. Probab.\/} 23(2):879--917.


\bibitem[Doeblin(1940)]{doeblin:1940} {\sc W.~Doeblin} (1940): Remarques sur la th\'{e}orie m\'{e}trique des fractions continues. \newblock {\em Compositio Mathematica\/} 7:353--371.

\bibitem[Freedman(1974)]{freedman:1974} {\sc D.~Freedman} (1974): The {P}oisson approximation for dependent events. \newblock {\em Ann. Probab.\/} 2:256--269.

\bibitem[Galambos(1972)]{galambos:1972}
{\sc J.~Galambos} (1972): The distribution of the largest coefficient in
  continued fraction expansions.
\newblock {\em Quart. J. Math.\/} 23(2):147--151.

\bibitem[Iosifescu(1977)]{iosifescu:1977}
{\sc M.~Iosifescu} (1977): A {P}oisson law for {$\psi $}-mixing sequences
  establishing the truth of a {D}oeblin's statement.
\newblock {\em Rev. Roumaine Math. Pures Appl.\/} 22:1441--1447.

\bibitem[Kallenberg(1983)]{kallenberg:1983}
{\sc O.~Kallenberg} (1983): {\em Random Measures\/}.
\newblock Akademie--Verlag, Berlin, 3rd edition.


\bibitem[Khintchine(1964)]{khintchine:1964} {\sc A.~Khintchine} (1964): Continued Fractions. \newblock{The University of Chicago Press}, Chicago, 1964.


\bibitem[Nakada and Natsui(2003)]{nakada:natsui:2003} {\sc H.~Nakada and R.~Natsui} (2003): On the Metrical Theory of Continued Fraction
Mixing Fibred Systems and Its Application to Jacobi-Perron Algorithm.
\newblock{Monatsh Math} 138: 267--288.


\bibitem[Philipp(1976)]{philipp:1976} {\sc W.~Philipp} (1976): A conjecture of Erd\"os on continued fractions.
\newblock{Acta Arithmetica} 28(4): 379--386.


\bibitem[Resnick(1987)]{resnick:1987}
{\sc S.~Resnick} (1987): {\em Extreme Values, Regular Variation and Point
  Processes\/}.
\newblock Springer-Verlag, New York.

\bibitem[Resnick(2007)]{resnick:2007}
{\sc S.~Resnick} (2007): {\em Heavy-Tail Phenomena \/}.
\newblock Springer-Verlag, New York.

\bibitem[Resnick and de Haan (1989)]{resnick:dehaan:1989}
{\sc S.~Resnick and L.~de Haan} (1989): {\em Second-order regular variation and rates of convergence in extreme-value theory \/}.
\newblock {\em Ann. Probab.\/} 24(1):97--124.

\bibitem[Smith(1988)]{smith:1988} {\sc R.~L.~Smith} (1988): Extreme value theory for dependent sequences via the Stein-Chen
method of Poisson approximation. \newblock {\em Stochastic Process. Appl.\/} 30(2):317--327. 

\bibitem[Tyran-Kami\'{n}ska(2010)]{tyran-kaminska:2010}
{\sc M.~Tyran-Kami\'{n}ska} (2010): Weak convergence to {L}\'{e}vy stable processes in dynamical systems.
\newblock {\em Stoch. Dyn.\/} 10(2):263--289.

\end{thebibliography}
\end{document}